\documentclass[11pt]{article}
\usepackage{}
\usepackage[total={6in, 9in}]{geometry}

\usepackage{amsfonts}
\usepackage{amsthm}
\usepackage{enumerate}
\usepackage{amssymb}
\usepackage{amsmath,amsfonts}
\usepackage{mathrsfs}
\usepackage{cases}
\usepackage{latexsym,bm}
\usepackage{indentfirst}
\usepackage{color}
\usepackage{ifpdf}
\usepackage{graphicx}
\usepackage{psfrag}
\usepackage{dsfont}
\usepackage[
pdfauthor={CESS},
pdftitle={Spanning trails},
pdfstartview=XYZ,
bookmarks=true,
colorlinks=true,
linkcolor=blue,
urlcolor=blue,
citecolor=blue,
bookmarks=true,
linktocpage=true,
hyperindex=true
]{hyperref}

\newtheorem{THM}{\textbf{Theorem}}

\newtheorem{CON}{\textbf{Conjecture}}
\newtheorem{COR}[THM]{\textbf{Corollary}}

\def\dist{{\fam0 dist}}

\begin{document}
\title{Towards the Small Quasi-kernel Conjecture}
\author{%
	Alexandr V. Kostochka  \thanks{Department of Mathematics, University of Illinois at Urbana-Champaign, Urbana, IL 61801, USA, and Sobolev Institute of Mathematics, Novosibirsk 630090, Russia. Email: {\tt kostochk@math.uiuc.edu}.  Supported by NSF grant DMS1600592, 
by Arnold O. Beckman Campus Research Board Award RB20003 of the University of Illinois at Urbana-Champaign and by grants 18-01-00353A and 19-01-00682 of the Russian Foundation for Basic Research.}
	%$^{\dagger}$%
	\and Ruth Luo \thanks{Department of Mathematics, University of California, San Diego, La Jolla, CA 92093 and University of Illinois at Urbana-Champaign,
	Urbana, IL 61801. Email: {\tt ruluo@ucsd.edu}. Supported by NSF grant DMS1902808.}
	\and Songling Shan 	\thanks{
	 Department of Mathematics, Illinois State University,
	Normal, IL 61790. Email: {\tt sshan12@ilstu.edu}.
}}
\date{\today}
\maketitle

 \begin{abstract} Let $D=(V,A)$ be a digraph. 
A vertex set $K\subseteq V$ is a {\em quasi-kernel} of $D$ if $K$ 
is an independent set in $D$ and for every vertex $v\in V\setminus K$, 
%there exist $w\in V\setminus K$ and $x\in K$ such that either $xv\in A$
%or $xw, wv\in A$, i.e., if
 $v$ is at most distance 2 from $K$. In 1974, Chv\'atal and Lov\'asz proved that every digraph 
has a quasi-kernel. P. L. Erd\H{o}s and L. A. Sz\'ekely in 1976 
conjectured that if every vertex of $D$ has a positive indegree, then $D$ 
has a quasi-kernel of size at most $|V|/2$.   This conjecture is only confirmed 
for narrow classes of digraphs, such as  semicomplete multipartite, quasi-transitive, or 
locally demicomplete digraphs. In this note, we state a similar conjecture for all digraphs, show that the two
conjectures are equivalent, and prove that both conjectures hold for a class of digraphs containing all
orientations of 4-colorable graphs (in particular, of all planar graphs).   

\end{abstract}
\vspace{2mm}

\section{Introduction  and notation}

%We consider digraphs whose underlying graphs are simple and finite. 
The digraphs in this note  may have antiparallel arcs, but do not have loops.
Let $D$ be a digraph.  We denote by $V(D)$ and $A(D)$ 
the vertex set and the arc set of $D$, respectively. 
We say $D$ is {\it weakly connected} if the underlying graph of $D$
is connected.
Let $x\in V(D)$.   
%Let $x,y\in V(D)$. If $xy\in A(D)$, we say $x$ \emph{dominates} $y$ and $y$ 
%is \emph{dominated} by $x$. 
The {\em open (closed) outneighborhood}  and {\em inneighborhood} of $x$ in $D$, denoted $N_D^+(x)$ ($N_D^+[x]$) and $N_D^-(x)$ ($N_D^-[x]$) are defined as follows. 
\begin{eqnarray*}
N_D^+(x)=\{y\in V(D)\,|\, xy\in A(D)\}, & N_D^+[x]=N_D^+(x)\cup\{x\}, \\
N_D^-(x)=\{y\in V(D)\,|\, yx\in A(D)\}, & N_D^-[x]=N_D^-(x)\cup\{x\}. 
\end{eqnarray*}
The {\em outdegree} of $x$ in $D$ is $d_D^+(x)=|N_D^+(x)|$, and the 
{\em indegree} of $x$ in $D$ is $d_D^-(x)=|N_D^-(x)|$.
Vertices of indegree zero in $D$ are called \emph{sources} of $D$
and vertices of outdegree zero in $D$ are called \emph{sinks} of $D$. 
By $\delta^+(D)$ (respectively,  $\delta^-(D)$) we denote the minimum outdegree (respectively, indegree) in $D$ among all vertices of $D$.
%and by $\delta^-(D)$ the minimum indegree in $D$ among all vertices of $D$. 
For each $X\subseteq V(D)$, we let 
\begin{eqnarray*}
	N_D^+(X)=\bigcup_{x\in X}N_G^+(x) \setminus X, & N_D^+[X]=N_D^+(X)\cup X,  \\
	N_D^-(X)=\bigcup_{x\in X}N_G^-(x) \setminus X, & N_D^-[X]=N_D^-(X)\cup X. 
\end{eqnarray*}

Let $u,v\in V(D)$ and $K\subseteq V(D)$. The distance from $u$ 
to $v$ in $D$, denoted $\dist_D(u,v)$, is the length of a shortest directed path from $u$
to $v$. 
 The distance from $K$ to $v$ 
in $D$,  is $\dist_D(K,v)=\min\{\dist_D(x,v)\,|\, x\in K\}$.  
We say $K$ is a \emph{kernel} of $D$ if $K$ is independent in $D$
and for every $v\in V(D)\setminus K$, $\dist_D(K,v)=1$. 
We say $K$ is a \emph{quasi-kernel} of $D$ if $K$ is independent in $D$
and for every $v\in V(D)\setminus K$, $\dist_D(K,v)\le 2$.

A digraph $D$ is \emph{kernel-perfect} if every induced subdigraph of  it  has a kernel. 
Richardson proved the following result.

\begin{THM}[Richardson~\cite{MR0014057}]\label{kernel}
	Every digraph without directed odd cycles is kernel-perfect. 
\end{THM}

  The proof gives rise to an algorithm to find one. On the other hand, Chv\'atal~\cite{kernelnp} showed that in general it is NP-complete to decide whether  a digraph has a kernel, and by a result of Fraenkel~\cite{MR675689} it is NP-complete even in the class of planar digraphs of degree at most 3.
While not every digraph has a kernel, % (for example, directed cycles of odd length), 
Chv\'atal and Lov\'asz~\cite{MR0414412} proved that every digraph 
has a quasi-kernel. In 1976, P.L. Erd\H{o}s and S. A. Sz\'ekely 
made the following conjecture on the size of a quasi-kernel in a digraph. 

\begin{CON}[Erd\H{o}s--Sz\'ekely~\cite{quasi-kernel}]\label{ESC}
	Every $n$-vertex digraph $D$ with $\delta^+(D)\ge 1$ has a quasi-kernel of size at most $\frac{n}{2}$. 
\end{CON}

If $D$ is an $n$-vertex digraph consisting of disjoint union of directed 2- and 4-cycles, then every kernel or quasi-kernel of $D$ has size exactly $\frac{n}{2}$.  Thus, Conjecture~\ref{ESC} is sharp.   

In 1996, Jacob and Meyniel~\cite{MR1395467} showed that a digraph without a kernel 
contains at least three distinct quasi-kernels. Gutin et al.~\cite{MR2051468}  characterized digraphs with exactly one and two-quasi-kernels, thus provided necessary and sufficient conditions for 
a digraph to have at least three quasi-kernels. However,  these results
do not discuss the sizes of the quasi-kernels. Heard and Huang~\cite{MR2419520} in 2008 showed that each digraph $D$ with $\delta^+(D)\ge 1$ has two disjoint quasi-kernels if $D$ is semicomplete multipartite (including tournaments), quasi-transitive (including transitive digraphs), or locally semicomplete.  As a consequence, Conjecture~\ref{ESC} is true for these three classes of digraphs.

We propose a  conjecture  which formally implies Conjecture~\ref{ESC}. It suggests a bound for  digraphs that may have sources.  Note that each quasi-kernel of a digraph 
 contains all of its source vertices and hence contains no outneighbors of the source vertices.
 
\begin{CON}\label{con2}
	Let $D$ be an $n$-vertex digraph, and let $S$ be the set of sources of $D$. Then $D$ has a quasi-kernel $K$ such that 
	$$
	|K|\le \frac{n+|S|-|N_D^+(S)|}{2}. 
	$$ 
\end{CON}

To show that the upper bound above is best possible, consider the following examples.
\begin{itemize}
\item Let $S$ be a nonempty set of 
isolated vertices, and let $D$ be a digraph obtained from a directed triangle by 
adding an arc from every vertex in $S$ to the same vertex in the triangle. Then  every  quasi-kernel of  $D$ has size $|S|+1=\frac{(|S|+3)+|S|-1}{2}$. 
\item Let $D$ be an orientation of a connected bipartite graph with parts $S$ and $T$ where each arc goes from $S$ to $T$. Then $S$  forms a quasi-kernel of $D$ of size $|S| = \frac{(|S|+|T|) + |S| - |T|}{2}$. 
\end{itemize}

In this paper, we support Conjectures~\ref{ESC} and~\ref{con2} by showing the following results.

\begin{THM}\label{4-partitea}
	Let $D$ be an $n$-vertex  digraph and $S$ be the set of sources of $D$. 
	Suppose that   $V(D)\setminus N_D^+[S]$
	has a partition $V_1\cup V_2$ such that $D[V_i]$
	is kernel-perfect
	 for each $i=1,2$. Then   $D$
	has a quasi-kernel of size at most $\frac{n+|S|-|N_D^+(S)|}{2}$.
\end{THM}

Since by  Theorem~\ref{kernel}, every  digraph without directed odd cycles is kernel-perfect, Theorem~\ref{4-partitea} immediately yields:

\begin{COR}\label{cor0}
 Conjectures~\ref{ESC} and~\ref{con2} hold for every orientation of each graph with chromatic number at most $4$.
\end{COR}

By the Four Color Theorem~\cite{4color1,4color2}, Corollary~\ref{cor0} yields that Conjectures~\ref{ESC} and \ref{con2}  hold for every digraph whose underlying graph is planar. 

\begin{THM}\label{THM:minc}
If  Conjecture~\ref{con2} fails and $D$ is a counterexample to it with the minimum number of vertices, 	then $D$ has no source.  
\end{THM}

%\begin{COR}\label{cor1}
%If  Conjecture~\ref{ESC} fails and $D$ is a counterexample to it such that   $D$ has the minimum number of vertices, 	then $D$  has a sink.  
%\end{COR}

Since Conjecture~\ref{con2} implies Conjecture \ref{ESC},  Theorem~\ref{THM:minc} %Corollary~\ref{cor1}
implies that the two conjectures are equivalent. 

%If a counterexample $D$ to Conjecture~\ref{con2} has a sink $v$ with at least two inneighbors, say $u_1,u_2$, then $D-u_1v$ also is a 
%counterexample  to Conjecture~\ref{con2}. Hence Theorem~\ref{THM:minc} immediately yields the following.
%
%\begin{COR}\label{cor2}
%If  Conjecture~\ref{con2} fails, then some counterexample $D$ to it with the minimum number of vertices
% has no source, has a sink, and every
% sink in $D$ has indegree~$1$.  
%\end{COR}

In the next section we prove Theorem~\ref{4-partitea} %which immediately yields Corollary~\ref{cor0} by Theorem~\ref{kernel}, 
and in Section~\ref{sec3} prove Theorem~\ref{THM:minc}.
%  and derive from it
%Corollary~\ref{cor1}.

\section{Proof of Theorem~\ref{4-partitea} }
%We will use the following result of Richardson~\cite{MR0014057}.
%\begin{LEM}[Richardson~\cite{MR0014057}]\label{kernel}
%	Every digraph without directed odd cycles has a kernel. 
%\end{LEM}

%\begin{thm1}\label{4-partite}
%	Let $D$ be an $n$-vertex  digraph and $S$ be the set of sources of $D$. 
%	Suppose that  $V(D)\setminus N_D^+[S]$
%	has a partition $V_1\cup V_2$ such that $D[V_i]$
%	has no directed odd cycles for each $i=1,2$. Then   $D$
%	has a quasi-kernel of size at most $\frac{n+|S|-|N_D^+(S)|}{2}$. 
%\end{thm1}

 % {\bf Proof of Theorem~\ref{4-partitea}.} 
Let $D_1=D-N_D^+[S]$, and 
 $V_1\cup V_2=V(D_1)$ be a partition of $V(D_1)$ such that $D[V_i]$  is kernel-perfect for each $i=1,2$.
In addition, we choose  such a partition so that $|V_2|$ is  as small as possible. Observe that adding a source vertex $v$ to a kernel-perfect digraph $H$
results in a new kernel-perfect digraph: let $H'$ be the resulting digraph, and let $F$ be a subdigraph of $H'$ that contains $v$. Then $K \cup \{v\}$ is a kernel of $F$ where $K$ is any kernel of $F - N^+_{H'}[v]$ in $H$. 

If there exists some $v \in V_2$ with no inneighbors in $V_1$, then we may move $v$ from $V_2$  to $V_1$, 
 and obtain a new partition of $V(D_1)$ into kernel-perfect subgraphs with a smaller $V_2$. % without creating any new odd cycles. 
Thus, by the choice of $V_2$, 
\begin{equation}\label{V2_indeDree}
N_{D_1}^-(v)\cap V_1\ne \emptyset \quad \mbox{ for every $v\in V_2$}. 
\end{equation}

Since $D[V_1]$  is kernel-perfect, it has a kernel $R$.  
Let  $R_0 $
be a smallest subset of $R$ such that $N_{D_1}^+(R_0)=N_{D_1}^+(R)$. 
By the minimality of $R_0$,  for every $v\in R_0$, there exists $u\in N_{D_1}^+(R)$
such that  $v$ is the only inneighbor of $u$ in $R_0$.
In particular, this yields
\begin{equation}\label{R0_size2}
|R_0|\le |N_{D_1}^+(R)|. 
\end{equation}
Let $D_2=D_1-(R_0\cup N_{D_1}^+(R))=D_1 - N_{D_1}^+[R_0]$, and   $S_2$ be the set of sources of $D_2$.  Since $R$ is a kernel in $D[V_1]$,
  $V(D_2) \subseteq V_2 \cup (R \setminus R_0)$. We partition $R \setminus R_0$ into the set $S_2'$ of the sources and 
  the set $T_2$ of the non-sources of $D_2$:  \[S_2':= S_2 \cap (R \setminus R_0), \qquad T_2=(R\setminus R_0) \setminus S_2.\]

%Since $V_1\setminus(R_0\cup N_{D_1}^+(R))=R\setminus R_0$, $T_2=(R\setminus R_0)\setminus S_2$. 
Then every vertex in $T_2$ has an inneighbor in the digraph $D_2$.  
As $V(D_2)=(R\setminus R_0)\cup (V_2\setminus N^+_{D_1}(R))$, and $T_2 \subseteq R$ is independent, we get 
\begin{equation}\label{T2}
N_{D_2}^-(v)\cap (V_2\setminus N_{D_1}^+(R))\ne \emptyset \quad \mbox{for every $v\in T_2$}. 
\end{equation}

%If $|T_2|\le |V(D_2)\setminus T_2|$, then let 
Let $
K=S\cup R_0\cup T_2=S\cup (R\setminus S_2).$
 Since $R_0\cup T_2\subseteq R$,  $K$ is an independent set in $D$. 
We now show that $\dist_D(K,v)\le 2$ for every $v\in V(D)\setminus K$. 
Note that 
$$ V_1\setminus K=S_2'\cup( N_{D_1}^+(R)\cap V_1), \quad \mbox{and}\quad V_2\cap K=\emptyset.$$

Since $K\supseteq S\cup R_0$, 
\begin{equation}\label{dist1}
\text{$\dist_D(K,v)=1$ for every $v\in N_D^+(S)\cup N_{D_1}^+(R)$. 
}
\end{equation}

Consider $v\in S_2'=R\cap S_2$. 
Since $v$ is not a source of $D$,  $N_D^-(v)\cap (N_D^+(S)\cup N_{D_1}^+(R))\ne \emptyset$. 
This, 
together with \eqref{dist1},  gives that $\dist_D(K,v)\le 2$. 
Lastly   let $v\in V_2\setminus N_{D_1}^+(R)$.
% By the choice of $R_0$,
%for every vertex $u\in S_2'$, 
%$N_{D_1}^+(u)\cap (V_2\setminus N_{D_1}^+(R))=\emptyset$.
%Because $N^+_{D_1}(R_0) = N^+_{D_1}(R)$, each vertex $u \in S_2'$ has no outneighbors in $V_2 \setminus N^+_{D_1}(R)$.  
By~\eqref{V2_indeDree}, $v$ must have an inneighbor in $V_1$. As $v \notin N^+_{D_1}(R)$, and every vertex in $V_1 \setminus R$ is an outneighbor of a vertex in $R_0$,
$\dist_D(K,v)= 2$. 

Therefore, $K$ is a quasi-kernel of $D$. If $|T_2|\le |V(D_2)\setminus T_2|$, then using~\eqref{R0_size2},
 $$|K|=|S|+|R_0|+|T_2|\leq |S|+\frac{1}{2}|R_0\cup N^+_{D_1}(R)|+\frac{1}{2}|T_2\cup V(D_2)|\leq \frac{1}{2}(n+|S|-|N_D^+(S)|),
 $$
 and the theorem holds.
Thus, assume that $|T_2|> |V(D_2)\setminus T_2|$ (so $|V(D_2)
\setminus T_2)| \leq |V(D_2)|/2$). 
Note that  $V(D_2)\setminus T_2=(V_2\setminus N_{D_1}^+(R))\cup S_2'$. 
Since 
  $D[V_2]$ is kernel-perfect and adding source vertices preserves kernel-perfectness,  
the digraph $D_2-T_2$ is also kernel-perfect. So let $W$ be a kernel of 
 $D_2-T_2$ and 
$K'=(S\cup R_0\cup W)\setminus N_D^+(W)$. 

Similarly to $K$, the set $K'$ is  independent in  $D$. Since   $|T_2|> |V(D_2)\setminus T_2|$,
$$|K'|\le |S|+|R_0|+|W|\leq |S|+\frac{1}{2}|R_0\cup N^+_{D_1}(R)|+\frac{1}{2}|V(D_2)|\leq
\frac{n+|S|-|N_D^+(S)|}{2}.$$
We now show that $\dist_D(K',v)\le 2$ for every $v\in V(D)\setminus K'$.  
Note that 
$$
S_2\subseteq K', \quad V_1\setminus K'= V_1 \cap [(R_0\setminus N_D^+(W))\cup (T_2\cup N_{D_1}^+(R))], \quad V_2\setminus K'=V_2\setminus W=N_{D}^+(W)\cap V_2. 
$$
By definition, $\dist_D(K',v)=1$ for every $v\in N_D^+(S)\cup N_{D_1}^+(R_0\setminus N_D^+(W)) \cup N_D^+(W)$.  
Thus, 
for every vertex $v\in N_{D_1}^+(R)\setminus N_{D_1}^+(R_0\setminus N_D^+(W))$, $\dist_D(K',v)\leq 2$. 
By~\eqref{T2}, for every $v\in T_2$, $\dist_D(K',v)\le 2$. 
Hence $K'$ is a quasi-kernel of $D$. 
\qed

\section{ Proof of Theorem~\ref{THM:minc}}\label{sec3} 
Assume
Conjecture~\ref{con2} fails and $D$ is a counterexample to it with the fewest vertices. % and modulo this, the fewest edges.
 Let $n=|V(D)|$. We assume  $n\ge 4$ as the cases $n \leq 3$ are verifiable by hand.  
By the minimality of $n$,  $D$ is weakly connected.  Let $S$ be the set of sources of $D$.  We show that $S=\emptyset$.  
 Assume instead that $S\ne \emptyset$.
 
 {\bf Case 1:}  $|N_D^+[S]|\ge 3$. 
 	 	Let $D_1$ be obtained from $D$ by  deleting all vertices in $N_D^+[S]$, 
 	adding two new vertices $x$ and 
 	$y$,  adding an arc from $y$ to every vertex of 
 	$D-N_D^+[S]$ that is an outneighbor of some vertex of $N_D^+(S)$ in $D$, and adding an arc from $x$ to $y$. 
 	Then $x$ is the only source vertex of $D_1$, and $N_{D_1}^+(x)=\{y\}$. 
 	Since $|V(D_1)|=|V(D)|-|N_D^+[S]|+2\le |V(D)|-1 $,
 	the minimality of $n$ implies that $D_1$
 	has a quasi-kernel $K_1$ of 
 	size at most $\frac{n-|N_D^+[S]|+2+1-1}{2}$. 
 	Then $K=(K_1\setminus \{x\})\cup S$ is a quasi-kernel of $G$
 	that has size at most 
 	$$
 	\frac{n-|N_D^+[S]|+2+1-1}{2}-1+|S|=\frac{n+|S|-|N_D^+(S)|}{2},
 	$$
 	as desired. 

 {\bf Case 2:}  $|N_D^+[S]|\ge 2$.  Since  $D$ is weakly connected, and $|S|\geq 1$, we get $|S|=1$ and $|N_D^+(S)|=1$.
 	 	Let $D_1=D-N_D^+[S]$. If $D_1$ has no sources, then  by the minimality of $D$, digraph $D_1$ has
  a quasi-kernel $K_1$  with 
 	$|K_1|\leq\frac{n-2}{2}$.  Then $K=K_1\cup S$ 
 	is a desired quasi-kernel of $D$. 
 	 	Therefore, we assume that $D_1$ has a source.
 	Let 
 	$$S_1=\{v\in V(D_1)\,|\, d_{D_1}^-(v)=0\}.$$
 	If $|N^+_{D_1}(S_1)|\le |S_1|$, we let  $D_2=D_1-S_1$.
 	By the minimality of $D$,  $D_2$
 	has a quasi-kernel $K_1$ of size at most
 	$\frac{n-2-|S_1|+|N_{D_1}(S_1)|}{2}\le \frac{n-2}{2}$.
 	Then $K=K_1\cup S$ 
 	is a desired quasi-kernel of $D$. 
 	Thus, we assume that $|N_{D_1}(S_1)|> |S_1|$. 
 	Let $D_2$ be obtained from $D_1$
 	by deleting all vertices in 
 	$N_{D_1}^+[S_1]$,  adding two new vertices $x$ and 
 	$y$,  adding an arc from $y$ to 
 	every vertex of 
 	$D_1-N_{D_1}^+[S_1]$ that is an outneighbor of some vertex of $N_{D_1}^+(S_1)$ in $D_1$, 
  and adding an arc from $x$ to $y$. 
 	Note that $x$ is the only source of $D_2$, and $N_{D_2}^+(x)=\{y\}$. 
 	Again, by the minimality of $D$,  $D_2$
 	has a quasi-kernel $K_1$ of 
 	size at most $\frac{n-2-|N_{D_1}^+[S_1]|+2+1-1}{2}$. 
 	Then $K=(K_1\setminus \{x\})\cup S\cup S_1$ is a quasi-kernel of $D$
 	that has size at most 
 	$$
 	\frac{n-2-|N_{D_1}^+[S_1]|+2+1-1}{2}-1+|S|+|S_1|\le \frac{n-1}{2},
 	$$
 	as desired. \qed

\bigskip
{\bf Acknowledgment.} We thank Peter L. Erd\H os for helpful discussions.

\end{document}